\documentclass{article}

\usepackage{PRIMEarxiv}

\usepackage[utf8]{inputenc} 
\usepackage[T1]{fontenc}    
\usepackage{hyperref}       
\usepackage{url}            
\usepackage{booktabs}       
\usepackage{amsmath,amsfonts, amsthm,amssymb}
\usepackage{algpseudocode}
\usepackage{algorithm}
\usepackage{bbm}
\usepackage{nicefrac}       
\usepackage{microtype}      
\usepackage{fancyhdr}       
\usepackage{enumitem}

\newcommand{\footnoteremember}[2]{
\footnote{#2}
\newcounter{#1}
\setcounter{#1}{\value{footnote}}
}
\newcommand{\footnoterecall}[1]{
\footnotemark[\value{#1}]
}

\newtheorem{theorem}{Theorem}
\newtheorem{definition}{Definition}
\newtheorem*{lemma}{Lemma}
\newcommand{\E}{\mathbb{E}}
\newcommand{\R}{\mathbb{R}}
\newcommand{\N}{\mathbb{N}}

\title{Computability of Optimizers
}

\author{
  Yunseok Lee\footnoteremember{LMU}{Mathematical Institute, Ludwig-Maximilians University of Munich, Munich, Germany} 
  \and Holger Boche\footnoteremember{TUM}{Institute of Theoretical Information Technology, Technical University of Munich, Munich, Germany}\footnoteremember{Bochum}{Excellence Cluster Cyber Security in the Age of Large-Scale Adversaries, Ruhr University Bochum, Bochum, Germany
  }\footnoteremember{MCQST}
  {Munich Center for Quantum Science and Technology (MCQST), Munich, Germany}
  \footnoteremember{MQV}{Munich Quantum Valley (MQV), Munich, Germany}
  \and Gitta Kutyniok\footnoterecall{LMU}\footnoteremember{Tromso}{Department of Physics and Technology, University of Troms\o, Troms\o, Norway}\footnoteremember{MCML}{Munich Center for Machine Learning (MCML), Munich, Germany}
}

\date{
\today}

\begin{document}
\maketitle

\begin{abstract}
Optimization problems are a staple of today’s scientific and technical landscape.
However, at present, solvers of such problems are almost exclusively run on digital hardware.
Using Turing machines as a mathematical model for any type of digital hardware, in this paper, we analyze fundamental limitations of this conceptual approach of
solving optimization problems.
Since in most applications, the optimizer itself is of significantly more interest than the optimal value of the corresponding function, we will focus on computability of the
optimizer.
In fact, we will show that in various situations the optimizer is unattainable on Turing machines and consequently on digital computers.
Moreover, even worse, there does not exist a Turing machine, which approximates the optimizer itself up to a certain constant error.
We prove such results for a variety of well-known problems from very different areas, including artificial intelligence, financial mathematics, and information theory, often
deriving the even stronger result that such problems are not Banach-Mazur computable, also not even in an approximate sense.

\end{abstract}


\section{Introduction}

Optimization is at the heart of basically any problem from science or industry; even the entire field of deep learning could not exist without optimization approaches. The numerical solvers of optimization problems are almost exclusively implemented on today's computers, i.e., on digital hardware such as CPUs or GPUs. However, the question of whether and which limitations this imposes is currently wide open. 
Since digital hardware is only able to handle discrete quantities (resp. bits) with arbitrary accuracy, solving optimization problems which admit continuous solutions requires approximations of the true solution. It is therefore of tremendous importance to mathematically analyze how large this misalignment between continuous-natured optimization problems and numerical solvers on digital hardware for those really is.  

In this paper we aim to provide a systematic approach to this problem from a computability viewpoint, studying whether there even do exist limits of computability due to the digital nature of the current hardware in contrast to the often continuous nature of optimization problems. Our results will unfortunately reveal that computability is indeed a major issue, as it will turn out that often the optimizers are {\em not computable} on digital hardware modeled by a Turing machine. 

\subsection{Optimization Problems} \label{subsec:opt}

Optimization is an area with a long and rich history. By optimization problems, we refer to the minimization or maximization of some functional $F: X \times Y \to \R$ over a solution space $X \subset \R^n$ and a parameter space $Y \subset \R^m$ for $n,m \in \N$, i.e.,
\begin{equation}\label{eq:opt}
\min_{x \in X(y)} F(x, y) \text{ or }\max_{x \in X(y)} F(x, y),
\end{equation}
where $y \in Y$ is a parameter in the parameter space and $X(y) \subset X$ is a subset of the solution space, depending on $y$. This general form allows to treat most problems from applications \cite{boyd}.

Two main problem settings can be identified in this context: The first asks for the optimal value or an approximation of it, i.e., constructing or approximating a function $\varphi : Y \to \R$ such that 
\begin{equation}\label{eq:optimalvalue}
\forall_{y\in Y}: \varphi(y) = \min_{x \in X(y)} F(x, y) \text{ or } \max_{x \in X(y)} F(x, y).
\end{equation}
The second problem setting aims to find an optimizer, i.e., to construct or approximate a function $G: Y \to X$ such that
\begin{equation}\label{eq:optimizer}
\forall_{y \in Y}: \varphi(y) = F(G(y), y),
\end{equation}
where $\varphi$ is the function defined by Equation \eqref{eq:optimalvalue}. It is evident that constructing a function $G$ yields a construction of a function $\varphi$. However, in general, the opposite direction does not hold. It is in this sense that finding $\varphi$ is ``easier'' than finding $G$.

Optimization problems suffer the same curse as many other problems of wide interest, namely there does in general not exist a closed-form solution for either $\varphi$ or $G$. Therefore, solutions usually have to be approximated by numerical algorithms run on today's computers. For a wide variety of optimization problems, 
established algorithms that aim to approximate the optimal solution do exist. A classical class of approaches are iterative solvers, which construct a sequence of approximators. In some cases, it has been proven that this sequence does indeed converge to the optimizer \cite{Arimoto}\cite{Portfolio}.

Depending on the application, either the function $G$ or the function $\varphi$ is of greater interest. Examples for the former are portfolio optimization or compressed sensing and, for the latter exemplary problems are computing the capacity of a channel or solving a deep learning problem. In practice, one often approximates $G$ by an iterative scheme to obtain a sequence $G_n : Y \to X$ and, correspondingly, $\varphi_n : Y \to \R$ through $\varphi_n(y):= F(G_n(y), y)$. Depending on the applied algorithm, one might obtain one or a combination of the following guarantees:
\begin{itemize}
    \item $\forall_{y\in Y} : \varphi_n(y) \to \varphi(y)$ with or without known convergence speed,
    \item $\forall_{y\in Y} : G_n(y) \to G(y)$ with or without known convergence speed.
\end{itemize}
Notice that in this case ``known convergence speed'' of a convergent Banach space sequence $a_n \to a$ refers to having an explicit description of a function $f: \N \to \R_+$ such that $\lim_{n\to \infty}f(n) = 0$ and $\|a_n - a\| \leq f(n)$ for all $n\in \N$.

\subsection{Computability}

Computability asks the question of how to mathematically model and analyze computations on perfect digital hardware. The term ``perfect'' refers to the assumption that there are no limitations regarding storage, computing power, and energy.
Additionally, perfect digital hardware is assumed to never make a mistake in a numerical calculation.  The only limitation is the fact that the number of calculation steps has to be finite.

It is evident that such a computability model is vastly superior to real-world digital hardware, i.e., real-world computers. Hence any practical algorithm can also be run on perfect digital hardware. Turing machines \cite{turing1936} are a version of perfect digital hardware and are considered the de-facto standard model for today’s digital computers. A Turing machine is a mathematical model of a machine capable of calculations by manipulating symbols on a strip of tape by only using a single reading head and a single internal state. Thus exploring the limitations of Turing machines allows to reveal the limitations of the computational abilities of today's most custom hardware. 

\subsection{Previous Work}

The earliest non-computability result has been provided by Church \cite{church} and Turing \cite{turing1936} \cite{turing1938}, by (independently of each other) proving non-computability of the Entscheidungsproblem. Since then the question of computability has been a staple in theoretical computer science and computer engineering and even found its way into the field of mathematical analysis \cite{weihrauch}.

Non-computability results for optimization problems have already been formulated in \cite{weihrauch}, \cite{boche2020algorithmic},  \cite{InverseProblems}, \cite{LinOpt},
\cite{noncomp1}, \cite{noncomp2}, \cite{InformationTheory}, and \cite{Colbrook2022}. 
The author of \cite{weihrauch} proves non-Borel-Turing computability for functions with certain discontinuities. Using these results \cite{LinOpt} deducts non-Borel-Turing computability of linear programs with real coefficients.
In \cite{InverseProblems} the stronger non-Banach-Mazur computability, as well as non-approximability, for inverse problems in the so-called lasso formulation is derived. The authors of \cite{InformationTheory}, \cite{noncomp1} and \cite{noncomp2} provide results on non-Banach-Mazur computability for selected problems in information theory and related fields.
Finally, \cite{Colbrook2022} proves non-computability for the specific case of finding neural networks to solve inverse problems. However, a more refined notion of non-computability is deployed, leading to a more nuanced result than the one obtained in this paper.

However, each of these results focuses on a special problem setting, leaving the question whether there exists a general comprehensive theory for non-computability for optimization problems wide open. In addition, most results are stated for non-Borel-Turing computability, yielding the question whether such a general theory can be formulated even for non-Banach-Mazur computability and non-approximability.

\subsection{Our Contributions}

In this paper, we develop such a general and comprehensive theory for non-computability or non-approximability for optimization problems in the more general setting of non-Banach-Mazur computability.

As mentioned in Subsection \ref{subsec:opt}, for a given optimization problem \eqref{eq:opt}, some approaches focus on finding the function $\varphi$ \eqref{eq:optimalvalue}, whereas others aim for computing $G$ \eqref{eq:optimizer}. Since \eqref{eq:optimizer} implies that computing $\varphi$ is ``easier'' than computing $G$, a significant amount of research in optimization has (successfully) focused on finding $\varphi$ over finding $G$. Naturally, the question arises whether $G$ can be non-computable even if $\varphi$ is computable. Hence the focus of this paper is on the question of computability of $G$ independently of $\varphi$.

As our main result, in Theorem \ref{mainthm}, we prove that surprisingly for a large class of optimization problems finding the optimizer, described by the function $G:Y \to X$, or even approximating it up to a constant error in a computable manner is not possible. While this is a result with far-reaching consequences (see Subsection \ref{subsec:impact}), it fortunately does not automatically imply that the optimal value, i.e., a function $\varphi: Y \to \R$ \eqref{eq:optimalvalue} or an approximation of it is non-computable. On the contrary, in Section 4 we present a selection of optimization problems, in which optimizers are non-computable and non-approximable, while the corresponding optimal values are in fact computable.

We also present and discuss several applications of Theorem \ref{mainthm}, showing non-computability and non-approximability for a selection of optimization problems (see Section 4). Intriguingly, most of such problems are even convex optimization problems. More detailed, we show non-Banach-Mazur computability and non-approximability for neural networks \ref{NN}, portfolio optimization \ref{Portfolio}, capacity maximizing distributions \ref{Information}, Wasserstein \ref{Wasserstein} distance, a lattice problem \ref{Lattice}, and linear programs \ref{LP}. Several of those examples also give evidence of the fact that our framework is comprehensive and includes some known results such as from \cite{LinOpt}, \cite{InformationTheory}, and \cite{Degen} as special cases.

\subsection{Impact of our Results} \label{subsec:impact}

We believe our results impact applications in several ways. In the following, we briefly discuss the most serious aspects.
\begin{itemize}
    \item \textit{Impossibility of general algorithms.} We prove non-existence of algorithms, which find or even approximate the optimizer, i.e., the search for such algorithms is in general futile. This even applies in the numerous cases, where the optimal value can be calculated or approximated by computable means. Especially for existing iterative algorithms, which are known to converge to the optimizer, our main result implies non-existence of a computable stop criterion, which ensures an arbitrarily small approximation error. 
    \item \textit{Limitations of digital hardware.} Our results are a consequence of considering general real-valued parameters instead of choosing a discrete space for $Y$. Consequently, non-computability results from the error caused by the digital nature of Turing machines aiming to approximate real numbers. This motivates the use of analog computer models such as the Blum-Shub-Smale machine \cite{BSS} or quantum computers \cite{quantumcomputers}.
    \item \textit{Importance of prior information.} Our results can be interpreted as a form of the ``no free lunch'' result \cite{NFL}. Without additional restrictions, for instance, concerning structure or regularity of the solution space $X$ or the parameter space $Y$, an algorithm, which finds or approximates an optimizer, cannot exist. This stresses the importance of prior knowledge.
\end{itemize}
We also believe that our results are of relevance for areas outside of optimization theory.
\begin{itemize}
\item \textit{Trustworthiness and Robustness Certificates.} In the context of trustworthiness of 6G-based communication \cite{trust6g1}\cite{trust6g2} and artificial intelligence \cite{trustNN}, our result is particularly concerning, since the absence of approximation guarantees might imply the absence of rigorous certificates. But the increased use of automatic systems makes trustworthy systems necessary, since large parts of sensitive areas, e.g., medicine, infrastructure, and autonomous vehicles and robots are hoping on these highly automated technologies.
\item \textit{Simulations.} Most simulations rely on forms of physics-based approximations. For those, it needs to be ensured that the approximated simulation is close to the real world. Our results imply that this type of approximation might not always be reliable. The negative impact of this problem is presumably even amplified in the context of synthetic data engines, e.g., Nvidia's Omniverse \cite{Omniverse}, which are simulations used to create data for machine learning tasks.
\end{itemize}

\subsection{Outline}

We introduce the computability framework in Section \ref{sec:comp}. 
In particular, we define Borel-Turing computability of a function and the more general Banach-Mazur computability notion of a function. We present our main theorem in Section \ref{sec:main} and provide conditions for non-Banach-Mazur and non-Borel-Turing computable functions in a general manner. In Section \ref{sec:App}, we then introduce a list of prominent optimization problems, with some background. All these problems will turn out to be non-Banach-Mazur computable and, in fact, not even approximable by a Banach-Mazur computable function.\\

\section{Introduction to Computability} \label{sec:comp}

In practice, most problems do not possess a closed-form solution.
Hence finding approximative solutions is a necessity.
Most approximative algorithms are designed to be run on Turing machines as an idealized model of today's digital hardware.
Historically, there has been a family of computing models, e.g., $\mu$-recursive functions \cite{mu} and $\lambda$-Calculus \cite{lambda}, which turned out to be equivalent to Turing machines. 
Although non-equivalent computing models exist, e.g., quantum computers\cite{quantumcomputers}, we will restrict ourselves to Turing machines, since digital hardware is the predominant hardware used today in real-life.\\
A comprehensive and formal introduction on the subject of computability can be found in \cite{pour-el}.
We start with a basic definition.
\begin{definition}
A function $f: \N \to \N$ is called \textbf{recursive} or \textbf{computable}, if there exists a Turing machine, which, given the input $x \in \N$, leaves $f(x)$ on its tape after termination.
With slight abuse of notation, we equate a recursive function with its corresponding Turing machine.
\end{definition}
\subsection{Computable Numbers}
Turing defined all rational numbers to be computable. The idea is that rational numbers can be used to approximate some real numbers arbitrarily well in a manner, which still allows digital computations. With this application in mind, Turing introduced the Turing machine in \cite{turing1936}.
We take over Turing's definition here.
\begin{definition}
A sequence of rational numbers $(r_n)_{n \in \N}$ is \textbf{computable}, if there exist recursive functions $s, p, q : \N \to \N$ such that
\[
\forall_{k \in \N}: q_k = (-1)^{s(k)} \frac{p(k)}{q(k)}.
\]
\end{definition}
Also, we adopt a version of convergence, which is more natural for Turing machines.
\begin{definition}
A sequence of real numbers $(x_n)_{n\in \N}$ does \textbf{converge effectively} to a limit $x \in \R$, if
\[
\forall_{n\in \N}: |x_n - x| \leq 2^{-n}.
\]
\end{definition}
\noindent Using the definition of computable rational sequences and effective convergence, we can define a real number to be computable if a Turing machine can approximate it with exponentially growing precision.
\begin{definition}
A real number $x \in \R$ is called \textbf{computable}, if there exists a rational computable sequence $(q_n)_{n \in \N}$, such that
\[
q_n \to x,
\]
where the convergence is effective.
The sequence $(q_n)_{n \in \N}$ is called a \textbf{representation} of $x$.
We refer to the set of computable real numbers $\R_c$.
\end{definition}
\noindent We remark that $\R_c \subset \R$ is a computable, dense set with field structure, i.e., closed under addition, subtraction, multiplication, and division, with exception of dividing through zero. Additionally, $\R_c$ is closed under effective convergence.\\
This definition is equivalent to the existence of a Turing machine, which outputs the base-2 representation of a real number up to the $n$th decimal place when getting $n$ on its input band.
Similarly, we define a real sequence to be computable if a Turing machine can approximate each member of the sequence with exponentially growing precision.
\begin{definition}
A sequence of real numbers $(x_n)_{n\in \N}$ is \textbf{computable}, if there exists a rational computable double sequence $(q_{n,k})_{n,k\in \N}$ such that
\[
\forall_{k,n \in \N}: | q_{n,k} - x_n | \leq 2^{-k}.
\]
\end{definition}
\noindent Notice that all these definitions can be extended to vector-valued quantities and sequences if every component satisfies the according definition.
\subsection{Computable Functions}
We start with the following definition of a computable function, which goes back to Turing \cite{turing1936} himself.
\begin{definition}
Let $N, M \in \N$.
A function $f: \R_c^N \to \R_c^M$ is \textbf{Borel-Turing computable}, if there exists a Turing machine, which transforms all representations $(r_n)_{n \in \N}$ of a vector $x \in \R_c^N$ to representations of $f(x)$.
\end{definition}
\noindent The following is a generalization of Borel-Turing computability, which involves the use of computable real sequences.
\begin{definition}
Let $N, M \in \N$.
A function $f:\R_c^N \to \R_c^M$ is \textbf{Banach-Mazur computable}, if for every computable real vector-valued sequence $(a_n)_{n\in \mathbb{N}}$, the sequence $(f(a_n))_{n \in \mathbb{N}} \subset \R^M$ is computable.
\end{definition}
\noindent We want to mention that all Borel-Turing computable functions are automatically also Banach-Mazur computable.
\subsection{Decidable Sets}
The most common definition of decidable sets concerns subsets of natural numbers.
\begin{definition}
A set $A \subset \N$ is called \textbf{decidable}, if the function $\mathbbm{1}_A: \N \to \N$, defined by
\[
\mathbbm{1}_A(n) := \begin{cases}
1, & n \in A,\\
0, & n \in A^c,
\end{cases}
\]
is recursive.
\end{definition}

\begin{definition}
A set $A \subset \N$ is called \textbf{semi-decidable}, if there exists a Turing machine $\mathcal{TM}_A$ such that $\mathcal{TM}_A(n)$ outputs 1, if $n \in A$, and $\mathcal{TM}_A(n)$ does not terminate, if $n \in A^c$.
\end{definition}
Note that the halting problem for Turing machines implies the existence of sets, which are semi-decidable but not decidable.
We can naturally extend this notion to subsets of $\R_c^n$ for any $n \in \N$.
\begin{definition}
Given $B\subset \R_c^n$, a set $A \subset B $ is called \textbf{(semi-)decidable}, if there exists a Turing machine $\mathcal{TM}_A$ such that $\mathcal{TM}_A(x)$ outputs 1, if $x \in A$ and $\mathcal{TM}_A(x)$ outputs 0 (resp. does not terminate), if $x \in B\backslash A$. $\mathcal{TM}_A(x)$ can either output an arbitrary symbol or not terminate for $x \in B^c$. 
\end{definition}
Intuitively, a Turing machine aims to decide if $x \in A$, using the prior information $x \in B$.
\section{Main Results} \label{sec:main}
\noindent
In this section, we develop a theory that allows to check for non-computability and even non-approximability of optimizers in a very flexible manner. 
Our results can be applied to a broad class of optimization problems with very different backgrounds, as we will see in Section \ref{sec:App}.\\
We now consider general optimization problems with some parameter space $Y \subset\mathbb{R}_c^m$ and solution space $X \subset \mathbb{R}_c^n$ over the continuous function $F: X \times Y \to \R_c$.
Additionally, with slight abuse of notation, we call $X(y) \subset X$ a subset of $X$ depending on $y \in Y$, i.e., we have a map $Y \to 2^X$, $y \mapsto X(y)$.
For example in the case of linear programs, $Y$ usually describes the linear inequalities, as well as the objective function and $X$, usually describes the space of solution vectors.\\
Given a fixed $y\in Y$, we are interested in optimization problems of the form 
\[
\min_{x \in X(y)}F(x,y) \text{ or } \max_{x \in X(y)}F(x,y).
\]
Most optimization problems can be written this way.\\
Usually, the existence of an optimizer is ensured through compactness of $X$ or $X(y)$. Note that compactness only proves abstract existence, but does not provide a description or approximation of the optimizer itself. 
In this paper, we use the following definition, if the optimization problem and its parameter space are clear,
\[
Opt(y) := \{x \in X(y) | F(x,y) =  \max_{x' \in X(y)}F(x',y)\}.
\]
The goal of optimization is to find a function $G: Y \to X$ such that 
\[
\forall_{y \in Y}: G(y) \in Opt(y)
\]
or at least an approximation of $G$, i.e., a function $G^* : Y \to X$ such that $G$ and $G^*$ are close. In our case we define closeness by
\[
\|G - G^* \|_\infty = \sup_{y \in Y} |G(y) -  G^*(y)| < \alpha
\] for some $\alpha > 0$. 
\begin{theorem}[Main Theorem]\label{mainthm}
Let $X, Y, X(y)$ and $F$ be as described above.
Let $G: Y \rightarrow X$ such that, for all $y \in Y$, we have $G(y) \in Opt(y)$.
Now let $Y_1, Y_2 \subset Y$,\: $y_1^* \in Y_1, y_2^* \in Y_2$ and $y_* \in Y$,  and $\gamma: [-1, 1] \rightarrow Y$, a Turing computable, continuous path such that:
\begin{enumerate}
    \item[(i)] $Y_1 \cap Y_2 = \varnothing$ and $G(Y_1) \cap G(Y_2) = \varnothing$,
    \item[(ii)] $\inf \limits_{y_1 \in Y_1\\ y_2 \in Y_2} \| y_1 - y_2 \| = 0$,
    \item[(iii)] $\inf \limits_{y_1 \in Y_1\\ y_2 \in Y_2} \| G(y_1) - G(y_2) \| = \kappa > 0$,
    \item[(iv)] $\gamma(-1) = y_1^*$, $\gamma(1) = y_2^*$, $\gamma(t_0) = y_*$,\\ for some $t_0 \in (-1, 1)$,
    \item[(v)] $\gamma([-1, t_0)) \subset Y_1$ and $\gamma((t_0, 1]) \subset Y_2$,
    \item[(vi)] $G(Y_1) \subset G(Y_1) \cup G(Y_2)$ is decidable.
\end{enumerate}
Then $G$ cannot be Borel-Turing computable.
In fact, there does not even exist a Borel-Turing computable function, which can approximate $G$ by up to an absolute error of $\alpha <  \frac{\kappa}{2}$, i.e. there does not exist a Borel-Turing computable function $G^*: Y \to X$ such that $\|G - G^* \|_\infty \leq \alpha$.
If we replace condition (vi) by \\
\\
\hspace{1cm} (vii) $Y_1 \subset Y_1 \cup Y_2$ is decidable,
\\
\\
then $G$ can even not be Banach-Mazur computable.
In fact, there does not even exist a Banach-Mazur computable function, which can approximate $G$ by up to an absolute error of $\alpha < \frac{\kappa}{2}$.
\end{theorem}
This has noteworthy consequences for computable stop criteria for iterative algorithms, which are guaranteed to converge to the optimizer.
Let $G: Y \to X$ and $\varphi: Y \to \R$ be defined as in \ref{eq:optimizer} and \ref{eq:optimalvalue}.
Assume for a given optimization problem that there exists an iterative scheme, that yields functions $G_n: Y \to X$ and consequently $\varphi_n: Y \to \R$, such that
\[
\forall_{y \in Y} : \varphi_n(y) \to \varphi(y) \text{ with known convergence speed and } \forall_{y \in Y} : G_n(y) \to G(y).
\]
This is the case for multiple optimization problems, examples being the Blahut-Arimoto algorithm for the capacity of a channel \cite{Blahut}\cite{Arimoto} and Cover's algorithm for portfolio optimization \cite{Portfolio}.
Now proving the non-computability and non-approximability of the function $G$ in this setup implies that even though there exist computable $G_n$, which converge pointwise to the sought function $G$, there can be \textit{no computable stop criterion}, which guarantees the error $\|G_n(y) - G(y)\|$ for all $y \in Y$ to be small.
Here being able to bound the error $\|\varphi_n(y) - \varphi(y)\|$ for all $y \in Y$ is irrelevant.
\section{Applications} \label{sec:App}
We can use Theorem \ref{mainthm} in a wide variety of cases.
In this section, we will provide a small sample of a few famous problems from a broad range of topics.
\subsection{Neural Networks}\label{NN}
Neural networks have seen a tremendous rise in popularity and successes in a wide variety of different areas, such as image processing \cite{image}, games \cite{atari}\cite{GO} and PDEs \cite{PDE}.
At the same time, there seem to be some inherent problems with neural networks, like instability \cite{instability}.
Consequently, the question if these problems might be an inherent property of neural networks was asked and some results seem to indicate that this is indeed the case \cite{Colbrook_2022}\cite{InverseProblems}.
We are continuing these results by showing a type of non-computability of neural networks, which, to our best knowledge, has been not shown before.
A \textit{(feed-forward) neural network} can be defined as functions $\Phi: \R^n \to \R^m$ of the form
\[
\Phi(x) := (A_L \circ \rho \circ A_{L-1} \circ \ldots \circ \rho \circ A_1)(x) 
\]
where $L \in \mathbb{N}$ and
\[
\forall_{l=1,\ldots ,L} : A_l x= W_lx + b_l, \: \: W_l \in \R_c^{n_l \times n_{l-1}}, b_l \in \R_c^{n_l}
\]
with $n_0 = n, n_L = m$ and $\rho : \R \to \R$ is a (non-linear) activation function, applied component-wise. 
For a more general theory on neural networks, we refer to \cite{ModernMathematics}.\\
In our setting we fix $\rho(x) = ReLU(x) = \max(0, x)$, for $x \in \R$, which is a very popular choice for neural networks \cite{jarrett}.
The matrix $A_l$ is called \textit{weight matrix} and $b_l$ is called \textit{bias vector}.
These are typically the free parameters of a neural network.
A choice of $L$ and $n_0,\ldots,n_L$ is called an \textit{architecture} of a neural network. 
We define the set of neural networks with architecture $n_0,\ldots,n_L$ as $\mathcal{NN}_{n_0,\ldots,n_L}$.
Now training a neural network with fixed architecture consists of minimizing a loss function over a data set $(x_i, y_i)_{i=1,\ldots,d}$, where $x_i \in \R^n$ and $y_i \in \R^m$ for $i=1, \ldots, d$.
A popular choice for a loss function to minimize is
\[
\min_{\Phi \in \mathcal{NN}_{n_0,\ldots,n_L}}W(\Phi, (x_i, y_i)_{i=1,\ldots,d}) : = \frac{1}{d}\sum_{i=1}^d \|\Phi(x_i) - y_i\|^2.
\]
Minimization is done only approximately by using a particular  form of gradient descent, namely, stochastic gradient descent \cite{BackProp}\cite{ADAM}.\\
For our setting, we ignore all biases, i.e. assume $b_l = 0$ for all $l$.
Additionally, we set all entries of the last layer, i.e., $A_L$, to 1.
For the following theorem, we consider the fixed architecture $L = 2,\: n_0 = 3,\: n_1 = 3,\: n_2 = 1$. 
\begin{theorem}[Neural Network] \label{thm:NeuralNetwork} 
Given $d \geq 14$ and a data set $\mathcal{D} = \prod_{i=1}^d(x_i, y_i) \in (\R^3_c \times \R_c)^d$, consider the minimization problem
\[
\min_{\Phi \in \R_c^{3\times 3} } \sum_{i=1}^d (\rho(\Phi x_i) - y_i)^2. 
\]
Let $G: (\R_c^{3} \times \R_c)^{d} \to \R_c^{3\times 3}$ such that, for all $\mathcal{D} \in (\R_c^{3} \times \R_c)^{d}$, we have: $G(\mathcal{D}) \in Opt(\mathcal{D})$.
Then G is not Banach-Mazur computable.\\
All functions $G^*:(\R_c^{3} \times \R_c)^{d} \to \R_c^{3\times 3}$ satisfying
\[
\|G - G^* \|_\infty \leq \alpha < 4
\]
are also not Banach-Mazur computable.
\end{theorem}

This theorem proves that no perfect loss-minimizing algorithms for neural networks in the special case of a shallow neural network with the architecture above can exist.
While this does not have to imply the same for wider and deeper neural networks, it is to be expected to also hold in more complicated cases.
Intuitively, calculating loss-minimizing neural networks gets "harder" with more parameters. Consequently, we expect similar results to hold true for general neural networks.
Note that this theorem does not imply non-approximability of neural networks interpreted as a function, but the non-approximability of its weights.
While this might be seen as a limitation of this theorem, since one is usually more interested in the neural network as a function itself instead of the precise weights, it cautions us against methods, which use or manipulate weights of a trained neural network directly as it is the case in, e.g., Dropout \cite{dropout} and Layer-Wise Relevance Propagation \cite{LRP}. Additionally, the non-approximability of weights points towards a fragile nature of neural networks and urges the need to ensure that the corresponding neural network function is indeed a good approximation of the desired function.
Also if there would be some computable way to always recover all possible weight configurations given a neural network function, this would imply non-approximability of the neural network function itself.
We believe even more non-computability results for deep learning exist and some inherent problems such as, e.g., instability \cite{instability} might be a fundamental flaw of neural networks on digital hardware, which are impossible to overcome completely.

\subsection{Financial Mathematics - Information Theory}\label{Portfolio}
In portfolio optimization, the stock market can be modeled by a random vector $X \in \R_+^m$, where each component describes a separate stock.
A portfolio $b\in \R_+^m$ is a vector such that $\sum_i b_i = 1$, which describes the allocation of the available funds.
The vector $b^tX$ describes the evolution of the portfolio after one time step.
Due to the multiplicative nature of investments, it is natural to maximize the convex functional \cite{Portfolio}\cite{latane}
\[
\max_{b}\E[\log(b^tX)],
\]
which is the expected return after one time step.
A well-known approach to this optimization problem is an iterative algorithm found by Cover \cite{Portfolio} in 1984. Cover's algorithm uses similar ideas as the Blahut-Arimoto algorithm \cite{Arimoto}\cite{Blahut} from information theory.\\
We will show that even though such an effective algorithm exists, finding a maximizing portfolio is non-computable in general. This implies no approximation guarantee for optimal portfolios in Cover's algorithm - or any other algorithm - can be made.\\
We consider the case of a discrete random vector $X(\cdot) = \sum_{i=1}^n \sum_{j=1}^m p_{i,j} \delta_{x_{i,j}}(\cdot) e_j$, where $p_{i,j} > 0$ are probabilities, i.e., $\sum_{i,j}p_{i,j} =1$, and $x_{i,j} \in \mathbb{R_+}$ are the possible outcomes.
Also $e_j \in \R^m$ are the standard basis vectors and we assume $\forall_{j,i_1 \neq i_2}: x_{i_1, j} \neq x_{i_2, j}$.
We define this set of discrete random vectors as $\mathcal{D}_{n,m}$.
\begin{theorem}[Log-Optimal Portfolio] \label{thm:portfolio}
Let $X \in \mathcal{D}_{n,m}$.
Define $W:\{b \in \mathbb{R}_{+}^m| \sum_i b_i = 1 \} \to \mathbb{R}$ by
\[
W(b) := \E [\log(b^t X)],
\]
and consider the corresponding maximization problem.
For all $n,m \in \mathbb{N}_+$ with $m>1$ define a function $G: \mathcal{D}_{n,m} \to \{b \in \mathbb{R}_{+}^m| \sum_i b_i = 1 \}$ such that for all  $y \in \mathcal{D}_{n,m}$, it holds $G(y) \in Opt(y)$.
Then $G$ is not Banach-Mazur computable.\\
All functions $G^*: \mathcal{D}_{n,m} \to \{b \in \mathbb{R}_{+}^m| \sum_i b_i = 1 \}$ satisfying
\[
\|G - G^* \|_\infty \leq \alpha < 1
\]
are also not Banach-Mazur computable.  
\end{theorem}
\subsection{Optimal Input Distribution - Information Theory}\label{Information}
In information theory, a point-to-point channel with one receiver and one transmitter is modeled by two discrete random variables $X$ and $Y$ over the probability spaces $\mathcal{X}$ and $\mathcal{Y}$. If we choose $\mathcal{X}$ and $\mathcal{Y}$ to be finite, we are describing a discrete memoryless channel (DMC). The channel itself is then given by a stochastic matrix $W \in \R^{m \times n}$, where $|\mathcal{X}| = n$ and $|\mathcal{Y}| = m$.
$X$, $Y$ and $W$ are related by 
\[
W(x) = P(Y|X=x).
\]
We define the mutual information of two discrete random variables $X$, $Y$ over $\mathcal{X}$, $\mathcal{Y}$ as 
\[
I(X,Y) = \sum_{x \in \mathcal{X}}\sum_{y \in \mathcal{Y}} P_{(X,Y)}(x,y) \log\left( \frac{P_{(X,Y)}(x,y)}{P_X(x)P_Y(y)} \right),
\]
where $P_{(X,Y)}$ is the probability mass function of $(X,Y)$, and $P_X$ and $P_Y$ are the probability mass functions of $X$ and $Y$.
Now the capacity $C(W)$ of a DMC $W$ is the maximal mutual information over all possible distributions over $\mathcal{X}$
\[
C(W) := \max_{X \in \mathcal{P}(\mathcal{X})} I(X,Y).
\]
The capacity of a DMC is well established and goes back to Shannon \cite{shannon}.
Also, more recently, the capacity has been considered in more complicated settings \cite{Boche1}\cite{Boche2}.
Trying to find the capacity of a DMC is a classical optimization problem for which a well-known approach using an iterative algorithm with convergence guarantee exists \cite{Arimoto}\cite{Blahut}.\\
We will show that even though such an effective algorithm exists, it is still impossible to compute a maximizing distribution in general or give an approximation guarantee.
This was already proven using the same construction in \cite{InformationTheory}, we will repeat the proof and show how this is a special case of Theorem \ref{mainthm}.
\begin{theorem}[Channel Capacity]\label{thm:channelcapacity}
Let $\mathcal{X}$ and $\mathcal{Y}$ be finite sets, such that $|\mathcal{X}| = n \geq 3$ and $| \mathcal{Y}| = m \geq 2$ and $W$ is a stochastic matrix. 
We define $\mathcal{P}(\mathcal{X}) := \{\text{random variables in }\mathcal{X} \}$ and $\mathcal{P}(\mathcal{Y}) := \{\text{random variables in }\mathcal{Y} \}$.
Since discrete random vectors can be identified by their probabilities for each event, i.e., we uniquely describe $Y \in \mathcal{P}(\mathcal{Y})$ by the vector $(P(Y=y_1), \ldots, P(Y=y_m)) \in \R_c^m$, with slight abuse of notation we equate those two objects by
\[
\mathcal{P}(\mathcal{Y}) \hat{=} \{ (v_1, \ldots, v_m) \in \R_c^m | \sum_{i=1}^m v_i = 1, \: \forall_i: v_i \geq 0 \}
\]
and analogously for $\mathcal{P}(\mathcal{X})$.
Also, define $\mathcal{W}$ as the set of all stochastic matrices in $\R_c^{m\times n}$.
Let $G: \mathcal{W} \to \mathcal{P}(\mathcal{X})$ be a function, such that, regarding the maximization problem
\[
\max_{X \in \mathcal{P}(\mathcal{X})} I(X, Y),
\]
for all $W \in \mathcal{W}$ we have $G(W) \in Opt(W)$.
Then G is not Banach-Mazur computable.\\
All functions $G^*: \mathcal{W} \to \mathcal{P}(\mathcal{X})$ satisfying
\[
\| G - G^* \|_\infty \leq \alpha < 1
\]
are also  not Banach-Mazur computable.
\end{theorem}
\subsection{Wasserstein Distance}\label{Wasserstein}
The Wasserstein-1 distance, originally formulated by Kantorovich \cite{Kantorovich} and Vaserstein \cite{Vaserstein} to tackle optimal transport problems, is a metric defined on the set of real probability distributions with finite first moment, i.e. ,
\[
\mathcal{P}_1 := \{\pi \text{ probability distributions}|\inf_{c\in \R} \int|x - c| d\pi(x) < \infty \}.
\]
One way to define the Wasserstein-1 metric $\mathbb{W}_1$ is by
\[
\mathbb{W}_1(\pi_1, \pi_2) := \sup_{f\in Lip_1} \left|\mathbb{E}_{x\in \pi_1}[f(x)] - \mathbb{E}_{x\in \pi_2}[f(x)]\right|,
\]
where $\pi_1, \pi_2 \in \mathcal{P}_1$ and
\[
Lip_1 := \{f:\mathbb{R} \to \mathbb{R} | \forall_{x,y \in \mathbb{R}}: |f(x) - f(y)| \leq |x-y| \}.
\]
This is the Kantorovich-Rubenstein duality formulation of the Wasserstein-1 distance.
Recently this formulation of the Wasserstein distance came to particular interest in the context of Wasserstein-GANs \cite{WGAN}.
The basic idea is to train a neural network, which is able to discriminate between the distribution of "nice" objects and the distribution of "adversarial" objects.
This is done by maximizing over $\left|\mathbb{E}_{x\in \pi_1}[f(x)] - \mathbb{E}_{x\in \pi_2}[f(x)]\right|$, where $f$ is the neural network to be trained and adding some regularizer to ensure that the Lipschitz constant is close to 1.
We consider the following relaxed setting.
First, we only consider probability distribution with computable density functions, supported in $[-\frac{1}{2}, \frac{1}{2}]$,
\[
P([-1/2, 1/2]) := \Big\{f:[-1/2, 1/2]\to \mathbb{R}_+ | \int f = 1,\; f\text{ Borel-Turing computable}\Big\}.
\]
Second, we restrict ourselves to a function space $\mathbb{F} \subset Lip_1$, which is made of Borel-Turing computable functions.
The only additional assumption on $\mathbb{F}$ is:
\[
\exists_{f_1, f_2\in \mathbb{F}}: \exists_{c_1, c_2\in \mathbb{R}}: \forall_{x\in [-1/2,1/2]} :f_1(x) = x + c_1 \text{ and }f_2(x) = |x| + c_2.
\]
This assumption holds in the example case of normalized neural networks.\\
We will show that calculating such a Wasserstein maximizer, or even approximating it is not possible with a Turing machine in these settings.
\begin{theorem}[Wasserstein distance]\label{thm:wasserstein}
We define
\[
\mathbb{W}'_1(\pi_1, \pi_2) := \sup_{f\in \mathbb{F}} \left|\mathbb{E}_{x\in \pi_1}[f(x)] - \mathbb{E}_{x\in \pi_2}[f(x)]\right|.
\]
Let $p_1, p_2 \in P[-\frac{1}{2}, \frac{1}{2}]$ be two computable probability densities.
Then the problem of finding a function $f \in \mathbb{F}$, such that $\mathbb{W}_1'(\pi_1, \pi_2) = \mathbb{E}_{x\in \pi_1}[f(x)] - \mathbb{E}_{x\in \pi_2}[f(x)]$ or $|\mathbb{W}_1'(\pi_1, \pi_2) - \mathbb{E}_{x\in \pi_1}[f(x)] + \mathbb{E}_{x\in \pi_2}[f(x)]| \leq \alpha < \frac{\sqrt{5}}{8\sqrt{3}}$ is not Banach-Mazur computable.
\end{theorem}
It might very well happen that such a maximizing function does not exist at all.
In this case, finding such a function is trivially non-computable.
However, we prove that finding such a maximizing function might not be computable even in the case a computable maximizing function does exist.
\subsection{Lattice Problem for Cryptographic Applications}\label{Lattice}
The basic idea of encryption is to apply a function $F$ to a message $m$ together with a (secret) key $k$ to obtain  an encrypted message $F(k,m) = e$.
Ideally, it is hard to recover $m$ from $e$ without knowledge of $k$ and easy to do with knowledge of $k$.
Usually, $F$ is motivated by using problems that are hard or suspected to be hard for all Turing machines to solve \cite{DiffieHellman} \cite{RSA} \cite{AES}.
So complexity and computability questions on Turing machines are central for well working encryption schemes.
We focus on the question of computability of one particular problem from the family of lattice problems.\\
Lattice problems have become a topic of interest with the rising feasibility of quantum computers.
Since quantum computers are able to crack conventional encryptions efficiently \cite{hallgren}\cite{shor}, lattice problems are seen as a new viable source for encryptions.
It is wildly believed, but not proven, that lattice problems are hard to solve not only for Turing machines \cite{LatticeNP} but also for quantum computers.\\
Different optimization problems are highly relevant candidates for post-quantum cryptography, among others \cite{regev} are the shortest vector problem, the shortest independent vector problem, the closest vector problem, and the short generator principal ideal problem.
Since Regev's discoveries \cite{regev}, tremendous efforts have been made to solve the mentioned problems \cite{eisentrager}\cite{biasse}. 
We consider the shortest independent vectors problem (SIVP), for which a randomized algorithm with exponential runtime exists if the complexity of the input is bounded \cite{sieve}. \\
To formulate the SIVP, we first have to define lattices over a field.
Commonly, finite fields such as $\mathbb{Z}/p\mathbb{Z}$, where $p\in \N$ is a large prime number, are considered for the message space.
In the following, we consider lattices over the field of real computable numbers $\R_c$.
As we will see, the transition from large $p$ to the "continuous" field $\R_c$ is problematic, and the corresponding optimization problem becomes non-computable on Turing machines, while for finite fields there have been recent successes \cite{eisentrager}\cite{biasse}.
Given $n \in \N$ and a basis $B = \{b_1 , \ldots, b_n \} \subset \R^n_c$, we define the corresponding lattice as 
\[
\Lambda(B) := \left\{\sum_{i=1}^n \lambda_i b_i \in \R_c^n| \forall_{i=1,\ldots, n}: \lambda_i \in \mathbb{Z}\right\}.
\]
Now define $\mathcal{B}(B)$ to be the set of bases in $\Lambda(B)$:
\[
\mathcal{B}(B) = \{ \beta \subset\Lambda(B) | \beta \text{ is basis of }\R_c^n\}.
\]
The SIVP  is described by the minimization problem
\[
\min_{\beta \in \mathcal{B}(B)} \sum_{b \in \beta} \|b\|_2.
\]
\begin{theorem}[SIVP]\label{thm:sivp}
Let $n\in \N$ and define $V$ as the set of all bases in $\R^n$.
Let $G: V \to V$ such that regarding SIVP and all bases $B\in V$ we have $ G(B) \in Opt(B)$.
Then G is not Banach-Mazur computable.\\
All functions  $G^* : V \to V$ satisfying 
\[
\| G - G^* \|_\infty \leq \alpha < \frac{\sqrt{2}}{2}
\]
are also not Banach-Mazur computable.
\end{theorem}
\subsection{Linear Program}\label{LP}
Linear programs are one of the most fundamental optimization problems in mathematics and are well-studied.
Although the most commonly used algorithm for linear programs is the simplex algorithm by Dantzig in 1947 \cite{simplex}, it is not the fastest algorithm in terms of  O-notation.
Indeed, Khachiyan proposed in 1979 an interior point algorithm \cite{LP}, which ensures polynomial runtime.
It is hard to understate the importance of this discovery for discrete mathematics. But some caveat remains until today --- most known algorithms deal only with integer coefficients.
And even in this case, the question remains if there exists an algorithm such that linear programs are solvable strongly polynomial.
The latter question has been cited by Smale in 1998 \cite{smale} among the 18 greatest problems for the 21st century.\\
For linear programs with real computable coefficients, we will show that the solution is not even Banach-Mazur computable.
We remark that the non-computability of linear programs has been already mentioned in \cite{LinOpt}.
\begin{theorem}[Linear program with real coefficients]\label{thm:lp}
Let $\mathcal{S}(A,y) := \{x \in \mathbb{R}_c^n | Ax \leq y \}$, where the inequality holds componentwise, with $A \in \mathbb{R}_c^{m\times n}, y \in \mathbb{R}_c^m$ and $c\in \mathbb{R}_c^n$.
Consider the maximization problem 
\[
\max_{x \in \mathcal{S}(A,y)} c\cdot x.
\]
Describing all coeffcients by the tuple $(A,y,c)$, we define the parameter space as $\mathbb{R}_c^{mn+m+n}$.\\
Let $G: \mathbb{R}_c^{mn+m+n} \to \mathbb{R}_c^n$ which outputs a maximizer for a given coefficient $(A, y, c)$, i.e. 
\[
\forall_{(A,y,c) \in \mathbb{R}_c^{mn+m+n}}: G(A,y,c)\in Opt(A,y,c).
\]
Then $G$ is not Banach-Mazur computable.\\
All functions $G^*: \mathbb{R}_c^{mn+m+n} \to \mathbb{R}_c^n $ satisfying
\[
\| G - G^* \|_\infty \leq \alpha < \frac{1}{2}
\]
are also not Banach-Mazur computable.
\end{theorem}

\section{Future Work}
We believe that our results can be extended in multiple directions.
\begin{itemize}
    \item \textit{Other optimization problems.} We believe more optimization problems are not computable.
    Using the characterization of Theorem \ref{mainthm}, it might be possible to characterize non-computable optimization problems even more precisely.
    \item \textit{Different computing models.} The question of computability for other computation models should allow to characterize the difference between Turing machines and other computation models more precisely. The question of stronger computation models could be vital to finding methods, which could actually calculate optimizers in a reliable and practically feasible way.
    
\end{itemize}

\section{Proof of Theorem \ref{mainthm}}\label{Proof of main theorem}
\begin{proof}
We start by assuming (vi) and showing that $G$ is not Borel-Turing computable.
Towards a contradiction - assume that there exists such a $G$, which is Borel-Turing computable.
WLOG we can assume $t_0 \neq \frac{k}{2^n}$ for all $n,k \in \mathbb{Z}$. Otherwise we rescale $\gamma$ slightly.\\
\\
Define $\mathcal{TM}_\gamma$ to be the Turing machine, which takes any representation $(r_1, r_2, \ldots)$ of a number $x \in \mathbb{R}_c\cap [-1,1]$ and outputs a representation of $\gamma(x)$.
We define $\mathcal{TM}_G$ analogously.\\
Then define $\mathcal{TM}^k_\gamma$ to be the Turing machine, which is identical to $\mathcal{TM}_\gamma$, but only requires the first $l$ numbers of a representation $(r_1, \ldots, r_l)$ of $x$ for a suitable large $l \in \mathbb{N}$.
Then it outputs the first $k$ numbers of a representation of $G(x)$ after finitely many steps.
If $\mathcal{TM}^k_\gamma$ does not get a suitable count of numbers of a representation to calculate $k$ numbers, it outputs an exception.
We define $\mathcal{TM}^k_G$ analogously.
\\
\\
Define $TM_{G(Y_1)}$ as the Turing machine, which takes the first $k$ numbers of a representation $(r_1,\ldots, r_k)$ of $x \in G(Y_1) \cup G(Y_2)$ and then outputs one of the following after a finite amount of steps:
\begin{enumerate}
    \item Decides that $x \in G(Y_1)$.
    \item Decides that $x \notin G(Y_1)$.
    \item Throws an exception, if $k$ is not large enough to decide.
\end{enumerate}
By assumption (vi), there exists some $k \in \mathbb{N}$ large enough such that $TM_{G(Y_1)}$ does not throw an exception.
This $k$ might depend on $x$.
\\
Now define a sequence $(a_n, b_n, c_n)_{n\in \N}$ recursively by
\begin{align*}
(a_1, b_1, c_1 ) &= (-1, 0, 1),\\
(a_{n+1}, c_{n+1}) &= \begin{cases}
(b_n, c_n), & \text{if } \gamma(b_n) \in Y_1\\
(a_n, b_n), & \text{if } \gamma(b_n) \in Y_2\\
\end{cases},\\
b_{n+1} &= \frac{a_{n+1} + b_{n+1}}{2}.
\end{align*}
\\
Then the following statements are true:
\begin{enumerate}
    \item $a_n$ is a computable rational sequence.
    \item $a_n \to t_0$ effectively.
    \item $\gamma(a_n) \in Y_1$ for infinitely many $n \in \N$.
    \item $\gamma(a_n) \in Y_2$ for infinitely many $n \in \N$.
\end{enumerate}
The same also holds for $b_n$ and $c_n$.
To prove statement 1 we define the following Turing machine:\\
\begin{algorithm}[H]
\caption{Recursively calculate $(a_n,b_n,c_n)_{n \in \mathbb{N}}$}
\begin{algorithmic}[1]
\Require $(a_n,b_n,c_n)$
\State $k \gets 1$
\State $l \gets 1$
\While{True}
\If{$\mathcal{TM}_G^k(\mathcal{TM}_\gamma^l(b_n))$ throws exception}
\State $l \gets l+1$
\Else
\State$(r_1,\ldots, r_k) \gets \mathcal{TM}_G^k(\mathcal{TM}_\gamma^l(b_n))$
\If{$TM_{G(Y_1)}(r_1, \ldots, r_k)$ throws exception}
\State $k \gets k+1$
\Else
\If{$(r_1, \ldots, r_k)$ is in $G(Y_1)$}
\State \Return $(b_n, \frac{b_n + c_n}{2}, c_n)$
\Else
\State \Return $(a_n, \frac{a_n + b_n}{2}, b_n)$
\EndIf
\EndIf
\EndIf

\EndWhile
\end{algorithmic}
\end{algorithm}
Note that $\mathcal{TM}_\gamma^l(b_n)$ never throws an exception, since $b_n$ is rational, so $b_n$ itself is a full representation of $b_n$.\\
\\

Statement 2 follows by using $|a_n - c_n| = 2^{-n+2}$, which implies $|t_0-a_n|, |t_0-c_n| \leq 2^{-n+2}$ and $|t_0 - b_n| \leq 2^{-n+1}$.\\
Statements 3 and 4 follow from our initial assumption $t_0 \neq \frac{k}{2^n}$ for any $k,n \in \mathbb{N}$.
So $t_0$ is computable with $(a_n)_{n\in \mathbb{N}}, (b_n)_{n\in \mathbb{N}}$ and  $(c_n)_{n\in \mathbb{N}}$ being representations.\\
This implies $y_*$ is computable using the representation $\mathcal{TM}_\gamma (b_1, b_2 , \ldots)$.
Define the subsequences $b_{\alpha_n}$ and $b_{\beta_n}$, such that
\begin{enumerate}
    \item $\{\alpha_n|n \in \mathbb{N}\} \cup \{\beta_n|n \in \mathbb{N}\} = \mathbb{N}$
    \item $\{\alpha_n|n \in \mathbb{N}\} \cap \{\beta_n|n \in \mathbb{N}\} = \varnothing$
    \item $(\alpha_n)_{n \in \mathbb{N}}, (\beta_n)_{n \in \mathbb{N}}$ are strictly monotonic increasing
    \item $\forall_n : \gamma(b_{\alpha_n}) \in Y_1 \wedge \gamma(b_{\beta_n}) \in Y_2$.
\end{enumerate}
By assumption (vi), the subsequences $b_{\alpha_n}$ and $b_{\beta_n}$ are computable, simply by adding a variable, which saves the condition checked in line 11 of algorithm 1. \\
Now by assumption (iii), either $\|G(y_*) - G(\gamma(b_{\alpha_n})) \| \geq \frac{\kappa}{2}$ or $\|G(y_*) - G(\gamma(b_{\beta_n})) \| \geq \frac{\kappa}{2}$ has to be true for all $n \in \N$ uniformly.\\
WLOG assume that for all natural $n$ it holds $\|G(y_*) - G(\gamma(b_{\alpha_n})) \| \geq \frac{\kappa}{2}$.
Let $A \subset \mathbb{N}$ be a semi-decidable but not decidable set with a matching Turing-machine $TM_A$.
Define the following computable double sequence:
\[
z_{n, k} := \begin{cases}
s,  &\text{if } TM_A(n) \text{ stops after } s \leq k \text{ steps} \\
k, &\text{else}
\end{cases}.
\]
Then consider the computable double sequence $(\tau_{n, k})_{n,k\in \mathbb{N}} := (b_{\alpha_{z_{n, k}}})_{n,k\in \mathbb{N}} \in Y_1$.
Define 
\[
\tau_n := 
\begin{cases}
\tau_{n, s}, & \text{if } TM_A(n) \text{ stops after } s \text{ steps}\\
t_0, &\text{else}
\end{cases}.
\]
Then we can show
\[
\tau_{n,k} \to \tau_n,
\]
where the convergence is effective.
This is obvious for $n \in A$ and for $n \notin A$ we can use $\|b_{\alpha_k}-t_0 \| \leq 2^{-k + 2}$.\\
So $(\tau_n)_{n \in \mathbb{N}}$ is a computable sequence.
Then $(G(\gamma(\tau_n)))_{n \in \mathbb{N}}$ is also a computable sequence by assumption and as an implication also $(\|G(\gamma(\tau_n)) - G(\gamma(t_0))\|)_{n \in \mathbb{N}}$.
Now notice that the following holds
\[
z_n := \|G(\gamma(\tau_n)) - G(\gamma(t_0))\|  \begin{cases}
\geq \frac{\kappa}{2}, &\text{if } n \in A\\
=0, &\text{if } n \notin A \\
\end{cases}.
\]
Using this property of $(z_n)_{n\in\N}$ we can define a Turing machine, which makes $A$ recursive.
Let $(r_{n,k})_{n,k \in \mathbb{N}}$ be a computable double sequence such that 
\[
|r_{n,k} - z_n| \leq 2^{-k}
\]
and $\mathcal{TM}_{r}$ be the Turing machine, for which $\forall_{n,k\in \mathbb{N}}: \mathcal{TM}_{r}(n,k) = r_{n,k}$.
Choose a computable $\delta \in \mathbb{R}_c$ s.t. $\frac{\kappa}{2}\geq\delta > 0$ and choose $B\in \mathbb{N}$ s.t. $2^{-B+1} < \delta$.
We call $TM_{< \delta}$ the Turing machine, which decides after a finite amount of steps, if an input is smaller than $\delta$. 
Otherwise it might run indefinitely.
Define $TM_{>\delta}$ analogously.
Now define the Turing machine
\begin{algorithm}
\caption{Check if $n \in A$ or $n \notin A$}
\begin{algorithmic}[1]
\Require $n \in \mathbb{N}$
\State $r_{n, B} \gets \mathcal{TM}(n,B)$
\State \textbf{run} $TM_{<\delta}(r_{n, B})$
\State \textbf{run} $TM_{>\delta}(r_{n, B})$
\If{$TM_{>\delta}(r_{n, B})$ terminates}
\State \Return $n \in A$
\EndIf
\If{$TM_{<\delta}(r_{n, B})$ terminates}
\State \Return $n \notin A$
\EndIf
\end{algorithmic}
\end{algorithm}
\\
This Turing machine decides for any $n\in\N$ if $n\in A$ or $n \notin A$. This implies the recursivity of $A$. So $G$ is not Turing-computable by contradiction.\\
\\
To prove non-approximability notice $\max (\|G(y_*) - G(\gamma(b_{\alpha_n})) \|,\|G(y_*) - G(\gamma(b_{\beta_n})) \|) \geq \frac{\kappa}{2}$.
Then consider the computable double sequence $(\tau_{n, k})_{n,k\in \mathbb{N}}$ defined by 
\[
\forall_{n,k \in \mathbb{N}} : \tau_{2n, k} := b_{\alpha_{z_{n, k}}} \in Y_1 \wedge \tau_{2n+1, k} := b_{\beta_{z_{n, k}}} \in Y_2
\]
Define 
\[
\tau_n = 
\begin{cases}
\tau_{n, s}, & \text{if } TM_A(\lfloor{\frac{n}{2}}\rfloor) \text{ stops after } s \text{ steps}\\
t_0, &\text{else}
\end{cases}
\]  
Then we can show
\[
\tau_{n,k} \to \tau_n,
\]
where the convergence is effective. 
So $(\tau_n)_{n \in \mathbb{N}}$ is a computable sequence.
Then $(G(\gamma(\tau_n)))_{n \in \mathbb{N}}$ is also a computable sequence, since $G$ is Turing-computable and as a consequence also Banach-Mazur computable, and as an implication also $(\|G(\gamma(\tau_n)) - G(\gamma(t_0))\|)_{n \in \mathbb{N}}$.
Now similarly to before we define:
\[
z_{n} = \|G(\gamma(\tau_{n})) - G(\gamma(t_0))\|.
\]
Then similarly to before it holds
\[
max(z_{2n}, z_{2n+1}) \begin{cases}
=0, &\text{if } n \in A \\
\geq \frac{\kappa}{2}, &\text{if } n \notin A
\end{cases}.
\]
Now assume there exists a Turing-computable function $G^*$, which approximates $G$ up to an absolute error $\alpha < \frac{\kappa}{2}$, i.e., $\|G^* - G\|_\infty \leq \alpha < \frac{\kappa}{2}$. Then, for $n \in A$ it holds
\begin{align*}
\kappa \leq &\| G(\gamma(\tau_{2n})) - G(\gamma(\tau_{2n+1}))\| \\
\leq & \|G(\gamma(\tau_{2n})) - G^*(\gamma(\tau_{2n}))\| \\
+ &\|G^*(\gamma(\tau_{2n})) - G^*(\gamma(\tau_{2n+1}))\| \\
+ &\|G^*(\gamma(\tau_{2n+1})) - G(\gamma(\tau_{2n+1})) \| \\
\leq & 2\alpha + \|G^*(\gamma(\tau_{2n})) - G^*(\gamma(\tau_{2n+1}))\|.
\end{align*}
Rearranging this inequality yields
\[
0 < \kappa - 2\alpha \leq \|G^*(\gamma(\tau_{2n})) - G^*(\gamma(\tau_{2n+1}))\|.
\]
Combining this result with the case $n \in A$ yields
\[
\|G^*(\gamma(\tau_{2n})) - G^*(\gamma(\tau_{2n+1}))\|
\begin{cases}
\geq \kappa - 2\alpha, &\text{if } n \in A \\
=0 , &\text{if } n \notin A
\end{cases}.
\]
Now choose a computable $\delta \in \mathbb{R}_c$ s.t. $0< \delta < \kappa- 2\alpha$, 
which makes $A$ recursive by using algorithm 2 as before.\\
\\
Now we prove the non-Banach-Mazur computability of $G$ under assumption (vii).
This case is similar to the last case but slightly simpler, since we don't have to use $\mathcal{TM}_G$ to calculate $b_n$.
Define the Turing machine
\begin{algorithm}[H]
\caption{Recursively calculate $(a_n,b_n,c_n)_{n \in \mathbb{N}}$}
\begin{algorithmic}[1]
\Require $(a_n,b_n,c_n)$
\State $k \gets 1$
\While{True}
\State$(r_1,\ldots, r_k) \gets \mathcal{TM}_\gamma^k(b_n)$
\If{$TM_{Y_1}(r_1, \ldots, r_k)$ throws exception}
\State $k \gets k+1$
\Else
\If{$(r_1, \ldots, r_k)$ is in $Y_1$}
\State \Return $(b_n, \frac{b_n + c_n}{2}, c_n)$
\Else
\State \Return $(a_n, \frac{a_n + b_n}{2}, b_n)$
\EndIf
\EndIf

\EndWhile
\end{algorithmic}
\end{algorithm}
Now we can repeat the same argumentation as before with algorithm 1, i.e., define computable sequences $(\tau_n)_{n\in \N}$, such that $\|G(\gamma(\tau_n)) - G(\gamma(t_0))\|$ is zero for $n \notin A$ and $\geq \frac{\kappa}{2}$ otherwise. Then using algorithm 2 one can prove the recursivity of $A$.
\end{proof}
\section{Proofs of applications}
\subsection{Neural Network}
\begin{proof}[Proof of Theorem \ref{thm:NeuralNetwork}]
The construction used in this proof is inspired by Example 2.5 in \cite{Degen}.

We use the term "neural network" and its weights in $\R^{3\times 3}$ interchangeably.
Define for $\epsilon \in [0, 1]\cap \R_c$:
\[
\gamma_1^\epsilon :=\left(1, 1, \frac{\epsilon}{2}\right),\; \gamma_2^\epsilon = \left(-1, 1, \frac{\epsilon}{3}\right),\; \gamma_3^\epsilon = \left(0, -2, \frac{\epsilon}{6}\right),
\]
\[
\Gamma_1^\epsilon := [\gamma_1^\epsilon| \gamma_2^\epsilon| \gamma_3^\epsilon]^T,
\]
\[
\Gamma_2^\epsilon := [-\gamma_1^\epsilon| -\gamma_2^\epsilon| -\gamma_3^\epsilon]^T.
\]

Now note that $\gamma_1^\epsilon + \gamma_2^\epsilon + \gamma_3^\epsilon = (0,0, \epsilon)$.
This implies for $\rho = ReLU$ and $x \in \mathbb{R}^3$,
\begin{align*}
\rho (\langle \gamma_1^\epsilon, x\rangle) +\rho (\langle \gamma_2^\epsilon, x \rangle) + \rho (\langle \gamma_3^\epsilon, x \rangle) = \rho (\langle -\gamma_1^\epsilon, x \rangle) + \rho (\langle -\gamma_2^\epsilon, x \rangle) + \rho (\langle -\gamma_3^\epsilon, x \rangle) + \langle (0, 0, \epsilon), x\rangle.
\end{align*}
We define the realization $\mathcal{R}: \R^{3\times 3} \to [\R^3 \to \R^3]$ of weights to be the corresponding neural network function, i.e.
\[
(\mathcal{R}([a|b|c]^T))(x) = \rho(\langle a, x\rangle) + \rho(\langle b, x\rangle) + \rho(\langle c, x\rangle).
\]
Define $F_1^\epsilon, F_2^\epsilon: \mathbb{R}^3 \to \mathbb{R}$ by $F_1^\epsilon = \mathcal{R}(\Gamma_1^\epsilon)$ and $F_2^\epsilon = \mathcal{R}(\Gamma_2^\epsilon) = \mathcal{R}(\Gamma_1^\epsilon) - \langle(0,0,\epsilon), \cdot\rangle$.\\
We introduce Lemma A.3 from \cite{Degen}, which in our special case can be written the following way.
\begin{lemma}[Lemma A.3 from \cite{Degen}]
Let $[a_1|a_2|a_3]^T,[a_1'|a_2'|a_3']^T \in \R^{3\times 3}$ be two neural networks with identical realization, i.e.,
\[
\forall_{x \in \R^3}: (\mathcal{R}([a_1|a_2|a_3]^T))(x) = (\mathcal{R}([a_1'|a_2'|a_3']^T))(x).
\]
If the vectors $a_1, a_2$ and $a_3$ are linearly independent, then there exist a permutation $\pi : \{1,2,3 \} \to \{ 1,2,3\}$, such that
\[
a_1 = a_{\pi(1)}', a_2 = a_{\pi(2)}', a_3 = a_{\pi(3)}'.
\]
\end{lemma}
Now since $\gamma_1^\epsilon, \gamma_2^\epsilon, \gamma_3^\epsilon$ are linearly independent, by Lemma A.3 from \cite{Degen} $\Gamma_1^\epsilon$ is the only way to parameterize the function $F_1^\epsilon$ with a neural network in $\R^{3\times 3}$ except permutations and trivial scalings of the last layer, which we can ignore since we set the last layer to be constantly 1.
The same holds for $\mathcal{R}(\Gamma_2^\epsilon)$.\\
\\
Now define the datasets
\begin{align*}
A_1^\epsilon = \Bigg\{ &\left(\left(\epsilon,0,-
4\right), 0\right), \left(\left(\epsilon,0,-
2\right),0\right), \left(\left(\epsilon,0,-1\right), 0.5\epsilon \right),\\
&((\epsilon,0,0), \epsilon), \left(\left(\epsilon,0,1\right), \frac{5}{3}\epsilon\right), \left(\left(\epsilon,0,3\right), 3\epsilon\right), \left(\left(\epsilon,0,6\right), 6\epsilon\right) \Bigg\},
\end{align*}
\begin{align*}
A_2^\epsilon = \Bigg\{ &\left(\left(\epsilon,0,-
4\right), 4\epsilon\right), \left(\left(\epsilon,0,-2\right),2\epsilon\right), \left(\left(\epsilon,0,-1\right), 1.5\epsilon \right),\\
&((\epsilon,0,0), \epsilon), \left(\left(\epsilon,0,1\right), \frac{2}{3}\epsilon\right), \left(\left(\epsilon,0,3\right), 0\right), \left(\left(\epsilon,0,6\right), 0\right) \Bigg\},
\end{align*}
\begin{align*}
B_1^\epsilon = \Bigg\{ &((0,\epsilon,-6), 0), ((0,\epsilon,-3),0), \left((0,\epsilon,-2.5),\frac{1}{6}\epsilon\right),\\
&\left((0,\epsilon,-2),\frac{1}{3}\epsilon\right), ((0,\epsilon,6), 7\epsilon), ((0,\epsilon,12), 12\epsilon), ((0,\epsilon,18), 18\epsilon) \Bigg\},\\
B_2^\epsilon = \Bigg\{ &((0,\epsilon,-6), 6\epsilon), ((0,\epsilon,-3),3\epsilon), \left((0,\epsilon,-2.5),\frac{8}{3}\epsilon\right),\\
&\left((0,\epsilon,-2),\frac{7}{3}\epsilon\right), ((0,\epsilon,6), \epsilon), ((0,\epsilon,12), 0), ((0,\epsilon,18), 0) \Bigg\},
\end{align*}
and $D_1^\epsilon = A_1^\epsilon \cup B_1^\epsilon$, as well as $D_2^\epsilon = A_2^\epsilon \cup B_2^\epsilon$.\\
We will show that the only neural networks in $\R^{3\times 3}$, which fit $D_1^ \epsilon$ resp. $D_2^\epsilon$ are $\Gamma_1^\epsilon$ resp. $\Gamma_2^\epsilon$, as well as trivial permutations of these.
We proof this for $D_1^\epsilon$, the analogous statement for $D_2^\epsilon$ follows by symmetry.\\
By Lemma A.3 from \cite{Degen} it suffices to proof that $\rho (\langle \gamma_1^\epsilon, x\rangle) +\rho (\langle \gamma_2^\epsilon, x \rangle) + \rho (\langle \gamma_3^\epsilon, x \rangle)$ is the only function realized by a neural network in $\R^{3\times3}$, which fits $D_1^\epsilon$.
So let $\Gamma \in \R_c^{3 \times 3}$ be a neural network s.t. $\forall_{(x,y) \in D_1^\epsilon }R(\Gamma)(x) = y$.
Now note, that the coordinates of $A_1^\epsilon$ all lie on the line $g = \{(1,0,t) \in \R^3 | t \in \R\}$.
Now $\mathcal{R}(\Gamma)|_g:g \to \R$ can be interpreted as a function $\R \to \R$, using the parametrization of $g$.
We write with slight abuse of notation for $t \in \R$
\begin{align*}
\mathcal{R}(\Gamma)|_g(t) = \mathcal{R}(\Gamma)((1,0,t)).
\end{align*}
Since $\Gamma$ has only 3 hidden neurons, $\mathcal{R}(\Gamma)|_g$ is a continuous piecewise linear function with at most 3 non-differentiable points.  
\\
Now by the choice of $A_1^\epsilon$, the only continuous piecewise linear function with at most 3 non-differentiable points fitting all data points of $A_1^\epsilon$ is $f: \R \to \R$ defined by
\[
f(t) := \begin{cases}
0 &, t \leq -\frac{2}{\epsilon}\\
1+t\frac{\epsilon}{2}&, -\frac{2}{\epsilon} < t \leq 0\\
1+t\frac{2\epsilon}{3} & ,0 < t \leq \frac{3}{\epsilon}\\
t\epsilon &, \frac{3}{\epsilon} < t
\end{cases}.
\]
This implies $f = \mathcal{R}(\Gamma)|_g$.
Now by using homogeneity of the ReLU function, i.e., for all $\lambda \in \R_+$ and $v \in \R^3$ we use $\mathcal{R}(\Gamma)(\lambda v) = \lambda\mathcal{R}(\Gamma)(v)$, we get for $x>0$ and $y\in \R$
\[
\mathcal{R}(\Gamma)(x,0,y) = \begin{cases}
0 &,y \leq -\frac{2}{\epsilon}x\\
x+\frac{\epsilon}{2}y &, -\frac{2}{\epsilon}x < y \leq 0\\
x +\frac{2\epsilon}{3}y&, 0 < y \leq \frac{3}{\epsilon}x\\
\epsilon y&, \frac{3}{\epsilon}x < y
\end{cases}. 
\]
Since $\Gamma$ has only 3 hidden nodes and the above description reveals the 3 lines of non-differentiability on the plane $E_1 = \{(x,0,y)| x,y\in \R  \}$
\begin{align*}
s_1 &= \{ \left(t\frac{\epsilon}{3},0,-t\right) | t\in\R\},\\
s_2 &= \{ \left(t\frac{\epsilon}{2},0,t\right) | t\in\R \},\\
s_3 &= \{ (t,0,0) | t \in\R\}.
\end{align*}
This implies
\[
\mathcal{R}(\Gamma)((x,0,y)) = \rho\left(\alpha \left(x - y\frac{\epsilon}{3}\right) \right) + \rho\left(\beta\left(x + y\frac{\epsilon}{2}\right)\right) + \rho(\gamma y)
\]
for some unknown $\alpha, \beta, \gamma \in \R$.\\
Since $\mathcal{R}(\Gamma)(x,0,y) = 0$ for all $y \leq -\frac{2}{\epsilon}x < 0$ this implies  $\alpha, -\beta, \gamma < 0$, i.e., in this region all hidden nodes are not activated.
Using $\mathcal{R}{\Gamma}(x,0,y) = x + y\frac{\epsilon}{2}$ for all $-\frac{2}{\epsilon}x < y \leq 0$, this implies in this region all but one hidden node are not activated since this is the region neighboring to the zero region above.
The activated neuron is clear by looking at the line of non-differentiability $\{ (x,y) \in \R^2 | y = -\frac{2}{\epsilon}x\}$. So we know that
\[
\mathcal{R}{\Gamma}(x,0,y) = \rho\left(\beta\left(x+y\frac{\epsilon}{2}\right)\right)
\]
for all $ -\frac{2}{\epsilon}x < y \leq 0$.
So this implies $\beta = 1$.
Repeating this line of argumentation we see that for $0 < y \leq \frac{3}{\epsilon} x$ now two neurons are activated.
The second neuron has to be the term $\rho(\gamma y)$ by looking at the line of non-differentiability again, i.e.,
\begin{align*}
\mathcal{R}{\Gamma}(x,0,y) = &\left(x+y\frac{\epsilon}{2}\right) + \frac{\epsilon}{6}y\\
= &\rho\left(\beta\left(x+y\frac{\epsilon}{2}\right)\right) + \rho(\gamma y),
\end{align*}
which implies $\gamma = \frac{\epsilon}{6}$.\\
Repeating this argument a final time shows $\alpha = -1$.
So it holds
\[
\mathcal{R}{\Gamma}(x,0,y) = \rho\left(x+y\frac{\epsilon}{2}\right) + \rho\left(y\frac{\epsilon}{3} - x\right) + \rho\left(\frac{\epsilon}{6}y\right).
\]
Now repeating the argument for $B_1^\epsilon $ yields
\begin{align*}
\mathcal{R}{\Gamma}(0,x,y) = \rho\left(x+y\frac{\epsilon}{2}\right) + \rho\left(x+y\frac{\epsilon}{3}\right) + \rho\left(-2x+\frac{\epsilon}{6}\right).
\end{align*}
Now since all coefficients of the z-coordinate (namely $\frac{\epsilon}{2},\frac{\epsilon}{3}$, and $\frac{\epsilon}{6}$) are unique there is a unique way to combine both descriptions to
\begin{align*}
\mathcal{R}{\Gamma}(x,y,z) = &\rho\left(x+y+z\frac{\epsilon}{2}\right) + \rho\left(-x+y+z\frac{\epsilon}{3}\right)\\
+ &\rho\left(-y+z\frac{\epsilon}{6}\right),
\end{align*}
which shows the uniqueness of the minimizer on $D_1^\epsilon$.\\
We can analogously proof the uniqueness of the optimizer for $D_2^\epsilon$.\\
\\
Now define $X = \R_c^{3 \times 3}, Y = (\mathbb{R}^3\times\mathbb{R})^d$ and $F: X \times Y \to \mathbb{R}$ by 
\[
F(\Gamma, ((x_1, y_1),\ldots, (x_d, y_d))) = \sum_{i=1}^d |R(\Gamma)(x_i) - y_i|^2.
\]
And define $Y_1 = \{D_1^\epsilon| 1>\epsilon > 0 \}$ and $Y_2 = \{D_2^\epsilon| 1>\epsilon > 0\}$ as well as the computable curve $\gamma: [-1, 1] \to Y$ by
\[
\gamma(x) = \begin{cases}
D_1^{(\frac{1}{2}-t)/100}, &\text{if } t \in [0, \frac{1}{2}]\\
D_2^{(t-\frac{1}{2})/100}, &\text{if }t \in (\frac{1}{2}, 1]
\end{cases}.
\]
We proved for $\epsilon > 0$  
\begin{align*}
Opt(D_1^\epsilon) = &\{ [\gamma_{\sigma(1)}^\epsilon |\gamma_{\sigma(2)}^\epsilon |\gamma_{\sigma(3)}^\epsilon  ]^T \in \R_c^{3\times 3} | \sigma \in S_3 \},\\
Opt(D_2^\epsilon) = &\{ [-\gamma_{\sigma(1)}^\epsilon |-\gamma_{\sigma(2)}^\epsilon |-\gamma_{\sigma(3)}^\epsilon  ]^T \in \R_c^{3\times 3} | \sigma \in S_3 \},
\end{align*}
where $S_3 = \{\sigma: \{1,2,3\} \to \{1,2,3\} | \sigma \text{ bijiective} \}$ is the permutation group of 3 elements.
So this implies
\[
\min_{y_1 \in Opt(D_1^\epsilon),y_2 \in  Opt(D_2^\epsilon)}\|y_1 - y_2 \|_1 = 8 > 0.
\]
Now applying Theorem \ref{mainthm} finishes our proof.
\end{proof}
\subsection{Financial Mathematics}
\begin{proof}[Proof of Theorem \ref{thm:portfolio}]
First note we can assume WLOG $m=2$ since for $m > 2$ we can choose the $x_{i,j} = \epsilon$ for arbitrary small $\epsilon >0$, $i = {1,\ldots, n}$ and $j > 2$.
In this case, the optimizer would ignore the returns for $X_j$, $j>2$.
Now choose for $\alpha > 0$, $x_{i,1} = i$ and $x_{i,2} = i\alpha$ as well as $p_{i,1} = p_{i,2} = \frac{1}{n}$.
Then it holds
\[
W(b) = \frac{1}{n}\sum_{i=1}^n \log \left(b_1 i + b_2 i\alpha\right).
\]
Now set $b_2 = 1-b_1$ and observe
\begin{align*}
W(b) =& \frac{1}{n}\sum_{i=1}^n \log \left(b_1i + (1-b_1)i\alpha \right)\\
=& \frac{1}{n}\sum_{i=1}^n \log \left(i(b_1 + (1-b_1)\alpha)\right).
\end{align*}
Here you can see
\[
Opt(X_\alpha) = \begin{cases}
\{(1,0)\},& \alpha < 1\\
\{(b_1, 1-b_1) | b_1 \in [0,1] \},& \alpha = 1\\
\{(0,1)\},& \alpha > 1
\end{cases}.
\]
So define $X = \{b \in \mathbb{R}_{+}^m| \sum_i b_i = 1 \}, Y = \{\text{discrete random vectors in }\mathbb{R}^m_+ \text{ with n outcomes}\}$, $Y_1 = \{X_\alpha | \alpha > 0\}$, and $Y_2  = \{X_\alpha | \alpha < 0\}$.
The computable path is defnied as $\gamma: [-1, 1] \to Y$, $\gamma(t) = X_t$.
$Y_1 \subset Y_1 \cup Y_2$ is decidable by checking if $\alpha$ is positive or negative. 
Now using 
\[
\inf_{y_1 \in Y_1, y_2 \in Y_2} \|Opt(y_1) - Opt(y_2) \|_1 = 2
\]
we can apply Theorem \ref{mainthm} to conclude the proof.
\end{proof}

\subsection{Information Theory}
\begin{proof}[Proof of Theorem \ref{thm:channelcapacity}]
We start with the case $n = 3$ and $m = 2$. Define
\[
W_* = \begin{pmatrix}
1 & 0 & 0\\
0 & 1 & 1
\end{pmatrix}
\]
and also
\[
W_{1, \mu} = \begin{pmatrix}
1 & 0 & \mu\\
0 & 1 & 1-\mu
\end{pmatrix},
\:
W_{2, \mu} = \begin{pmatrix}
1 & \mu & 0\\
0 & 1 - \mu & 1
\end{pmatrix}
\]
for $\mu \in (0, 1)$.
Let
\[
\mathcal{P}_1 = \left\{ (p_1,p_2,p_3) \in \mathcal{P}(\mathcal{X}) | p_1 = \frac{1}{2}, \: p_2 + p_3 = \frac{1}{2} \right\}
\]
and
\[
\mathcal{P}_1 = \left\{ (p_1,p_2,p_3) \in \mathcal{P}(\mathcal{X}) | p_2 = \frac{1}{2}, \: p_1 + p_3 = \frac{1}{2} \right\}.
\]
We now show $Opt(W_{1,\mu}) = \mathcal{P}_1$ and $Opt(W_{2,\mu}) = \mathcal{P}_2$.
For this we define 
\[
W_* = \begin{pmatrix}
1 & 0 & 0 \\
0 & 1 & 1 
\end{pmatrix}, \;
\hat{W} = \begin{pmatrix}
1 & 0 & 1 \\
0 & 1 & 0
\end{pmatrix}.
\]
For $p=(p_1, p_2, p_3) \in \mathcal{P}(\mathcal{X})$ we consider
\begin{align*}
I(p, W_*) = &p_1 \cdot 1 \cdot \log \frac{1 \cdot p_1}{p_1 \cdot p_1} + p_2 \cdot 1 \cdot \log \frac{1 \cdot p_2}{p_2 (p_2 + p_3)}\\
 + &p_3 \cdot 1 \cdot \log \frac{1 \cdot p_3}{p_3(p_2 + p_3)}\\
=&p_1 \cdot 1 \cdot \log \frac{1 \cdot p_1}{p_1 \cdot p_1} + p_2 \cdot 1 \cdot \log \frac{1 \cdot p_2}{p_2 (p_2 + p_3)}\\
 + &p_3 \cdot 1 \cdot \log \frac{1 \cdot p_3}{p_3(p_2 + p_3)}\\
= &p_1 \log\frac{1}{p_1} + (p_2 + p_3)\log \frac{1}{p_2 + p_3} \\
= & p_1 \log\frac{1}{p_1} + (1 - p_1) \log \frac{1}{1-p_1} \\
= & h_2(p_1)
\end{align*}
Where $h_2$ is the binary entropy function.
Is is well known that $h_2(x)$ is maximal if and only if $x= \frac{1}{2}$. 
This immediately implies
\[
Opt(W^*) = \mathcal{P}_1
\] 
and analogously 
\[
Opt(\hat{W}) = \mathcal{P}_2.
\]
Now using 
\[
W_{1, \mu} = (1 - \mu)W_* + \mu \hat{W}
\]
as well as convexity it holds for $p \in \mathcal{P}(\mathcal{X})$
\begin{align*}
I(p, W_{1, \mu}) \leq (1 - \mu) I(p, W_*) + \mu I(p, \hat{W}).
\end{align*}
So for $\mu \in (0, 1)$ we have
\[
Opt(W_{1, \mu}) = \mathcal{P}_1 \cap \mathcal{P}_2 = \left\{\begin{pmatrix}
\frac{1}{2} \\ \frac{1}{2} \\ 0
\end{pmatrix}
\right\}.
\]
Analogously we can proof for $\mu \in (0, 1)$
\[
Opt(W_{2, \mu}) = \left\{\begin{pmatrix}
\frac{1}{2} \\ 0 \\ \frac{1}{2} 
\end{pmatrix}
\right\}.
\]
Now we can set $X = \mathcal{P}(\mathcal{X})$, $Y = \{M\in \R_c^{m\times n} | M \text{ is stochastic matrix}\}$, $Y_1 = \{W_{1, \mu} | \mu \in (0,1) \}$, $Y_2 = \{W_{2, \mu} | \mu \in (0,1) \}$ as well as $\gamma : (-1, 1) \to Y$ defined by
\[
\gamma(t) = \begin{cases}
W_{1, -t}, & t < 0\\
W_*, & t = 0\\
W_{2, t}, & t > 0
\end{cases}.
\]
It holds
\[
\inf_{y_1 \in Opt(Y_1), y_2 \in Opt(Y_2)} \| y_1 - y_2 \|_2 > 2.
\]
Since $Y_1 \subset Y_1 \cup Y_2$ is decidable by checking if for $W \in Y_1 \cup Y_2$ either $(W)_{1,2}>0$ or $(W)_{1,3}>0$, theorem \ref{mainthm} finishes the proof for $n=3$ and $m=2$.\\
For the cases of $n >3$ or $m>2$, we can reduce them to the case we just proved by using the following construction.
Let $\overline{W} \in \R^{2 \times 3}$ be a stochastic matrix for a DMC.
Let $W \in \R^{m \times n}$ be the a stochastic matrix defined by
\[
W(y|x) = \begin{cases}
\overline{W}(y|x), & y \in \{1, 2\}, x \in \{1, 2, 3\}\\
\overline{W}(y|1), & y \in \{1,\ldots, m \}, x \in \{4, \ldots, n \}\\
0, & y \in \{3, \ldots, m \}, x \in \{1 , \ldots, n \}\\
\end{cases}
\]
Now assume $G: \mathcal{W} \to \mathcal{P}(\mathcal{X})$ to be a function, such that for all $W \in\mathcal{W}$ we have $G(W) \in Opt(W)$.
Now let 
\[
G(W) = \begin{pmatrix}
g_1(W)\\
\vdots\\
g_n(W)
\end{pmatrix}
\]
and define the functions 
\begin{align*}
g_1^*(\overline{W}) &:= g_1(W) + \sum_{x=4}^n p_x(W),\\
g_2^*(\overline{W}) &:= g_2(W),\\
g_3^*(\overline{W}) &:= g_3(W).\\
\end{align*}
Now consider the function 
\[
G'(\overline{W}) :=\begin{pmatrix}
g_1^*(\overline{W})\\
g_2^*(\overline{W})\\
g_3^*(\overline{W})
\end{pmatrix},
\]
which is a composition of the following components:
1. the function $\overline{W} \to W$, described by the construction above;
2. the function $G$;
3. the functions $g_i^*$ for $i=1,2,3$.\\
Note that the first and third component is Banach-Mazur computable. So $G'$ is Banach-Mazur computable if and only if the function $G$ is Banach-Mazur computable.
But we have $G'(\overline{W}) \in Opt(\overline{W})$, which implies that $G'$ cannot be Banach-Mazur computable and consequently $G$ cannot be Banach-Mazur computable.\\
Analogously we can show, that any function $G^*: \mathcal{W} \to \mathcal{P}(\mathcal{X})$, satisfying $\|G - G^* \|_\infty < 1$ cannot be Banach-Mazur computable.
\end{proof}
\subsection{Wasserstein Distance}
\begin{proof}[Proof of Theorem \ref{thm:wasserstein}]
Define for probabilty distributions $\pi_1$ and $\pi_2 \in P[-1/2, 1/2]$ and $f \in \mathbb{F}$ the function 
\[
W(\pi_1, \pi_2, f) := |\mathbb{E}_{x\in \pi_1}[f(x)] - \mathbb{E}_{x\in \pi_2}[f(x)]|.
\]
Define the following equivalence relation $\sim$ on $\mathbb{F}$:\\
For $f, g\in \mathbb{F}$, $f\sim g$ holds iff
\[
\exists_{c\in \mathbb{R}}: f = g + c \; \lor f = -g + c.
\]
Define the corresponding equivalence class $[f] = \{ g \in \mathbb{F} | f \sim g\}$ and $\mathbb{F}/\sim := \{[f] | f\in \mathbb{F} \}$.\\
Now note that $f \sim g$ implies $W(\pi_1, \pi_2, f) = W(\pi_1, \pi_2, g)$.\\
Define for $\epsilon \geq 0$ the following density functions, which are in $P[-1/2,1/2]$:
\[
p_{1,\epsilon}(x) =
1 + \epsilon x,
\]
as well as 
\[
p_{2,\epsilon}(x) = \begin{cases}
1+\epsilon+4\epsilon x, & -\frac{1}{2} \leq x < 0\\
1+\epsilon-4\epsilon x, &  0\leq x \leq \frac{1}{2}
\end{cases},
\]
and 
\[
p_*(x) = p_{1,0}(x) = p_{2,0}(x) = \mathbbm{1}_{[-\frac{1}{2}, \frac{1}{2}]}(x),
\]
and $\pi_{1,\epsilon}, \pi_{2,\epsilon}, \pi_*$ to be the corresponding probability distributions.
Now notice that 
\[
Opt(\pi_1,\pi_2) = \{f\in \mathbb{F}| W(\pi_1, \pi_2, f) = \mathbb{W}'_1(\pi_1, \pi_2) \}.
\]
We want to show
\[
Opt(\pi_{1,\epsilon}, \pi_*) = [id|_{[-\frac{1}{2}, \frac{1}{2}]}]
\]
and
\[
Opt(\pi_{2,\epsilon}, \pi_*) = [|\cdot||_{[-\frac{1}{2}, \frac{1}{2}]}].
\]
\textbf{Proof of $Opt(\pi_{1,\epsilon}, \pi_*) = [id|_{[-\frac{1}{2}, \frac{1}{2}]}]$:}\\
For any $f \in \mathbb{F}$ it holds
\begin{align*}
&W(\pi_{1,\epsilon}, \pi_*, f) = \epsilon\left|\int_{-\frac{1}{2}}^{\frac{1}{2}} f(x) x dx \right|\\
=& \epsilon\left|\int_{-\frac{1}{2}}^{\frac{1}{2}} (f(x) - f(0)) x dx \right|\\
\leq& \epsilon \int_{-\frac{1}{2}}^{\frac{1}{2}} x^2\\
=& \frac{\epsilon}{12}.
\end{align*}
Notice that $W(\pi_{1,\epsilon}, \pi_*, f) = \frac{\epsilon}{12}$ if $f \in \mathbb{F}$ with $f(x) = x + c$ for all $x \in [-\frac{1}{2}, \frac{1}{2}]$ and some $c \in \R$.
By assumption such a function $f\in \mathbb{F}$ exists.
So the inequality above becomes an equality for all maximizers $f \in \mathbb{F}$.
This implies that the only maximizers are functions with constant slope $\pm 1$ in the domain $[-\frac{1}{2}, \frac{1}{2}]$, i.e., all functions in $[id_{[-1/2, 1/2]}]$.\\
\textbf{Proof of $Opt(\pi_{2,\epsilon}, \pi_*) = [|\cdot||_{[-\frac{1}{2}, \frac{1}{2}]}]$:}\\
Again we calculate for arbitrary $f \in \mathbb{F}$:
\begingroup
\allowdisplaybreaks
\begin{align*}
&W(\pi_{2,\epsilon}, \pi_*, f) = \left| \int_{-\frac{1}{2}}^\frac{1}{2} f(x)(1-p_{2,\epsilon}(x)) dx \right|\\
=&\epsilon\left|\int_{-\frac{1}{2}}^{0} f(x)(1 + 4x) dx + \int_{0}^{\frac{1}{2}} f(x)(1 - 4 x) dx \right|\\
=&\epsilon\Bigg|\int_{-\frac{1}{2}}^{0} \left(f(x) - f(-1/4)\right)(1 + 4x) dx + \int_{0}^{\frac{1}{2}} (f(x) - f(1/4))(1 - 4 x) dx \Bigg|\\
\leq& \epsilon \int_{-\frac{1}{2}}^{0} |f(x) - f(-1/4)||1 + 4x| dx + \epsilon\int_{0}^{\frac{1}{2}} |f(x) - f(1/4)||1 - 4 x| dx \\
=& \epsilon\Bigg( \int_{-\frac{1}{2}}^{-\frac{1}{4}} |f(x) - f(-1/4)|(-1 -  4x) dx + \int_{-\frac{1}{4}}^{0} |f(x)-f(-1/4)|(1 + 4x)  dx\\
+& \int_{0}^{\frac{1}{4}} |f(x)-f(1/4)|(1 - 4x)  dx + \int_{\frac{1}{4}}^{\frac{1}{2}} |f(x) - f(1/4)|(4x - 1)  dx \Bigg)\\
\leq& \epsilon\Bigg( \int_{-\frac{1}{2}}^{-\frac{1}{4}} |x + 1/4|(-1 -  4x) dx + \int_{-\frac{1}{4}}^{0} |x + 1/4 |(1 + 4x)  dx\\
+& \int_{0}^{\frac{1}{4}} |x - 1/4|(1 - 4x)  dx + \int_{\frac{1}{4}}^{\frac{1}{2}} |x - 1/4|(4x - 1)  dx \Bigg)\\
\leq& \epsilon\Bigg( \int_{-\frac{1}{2}}^{-\frac{1}{4}} (x + 1/4)(1 +  4x) dx + \int_{-\frac{1}{4}}^{0} (x + 1/4 )(1 + 4x)  dx\\
+& \int_{0}^{\frac{1}{4}} (x - 1/4)(4x - 1)  dx + \int_{\frac{1}{4}}^{\frac{1}{2}} (x - 1/4)(4x - 1)  dx \Bigg)\\
=& \epsilon\left( \int_{-\frac{1}{2}}^{0} (x + 1/4)(1 +  4x) dx + \int_{0}^{\frac{1}{2}} (x - 1/4)(4x - 1)  dx \right)\\
=& \frac{\epsilon}{12}.
\end{align*}
\endgroup
Notice that $W(\pi_{2, \epsilon}, \pi_*, f) = \frac{\epsilon}{12}$ for $f \in \mathbb{F}$ s.t. $\forall_{x\in [-\frac{1}{2}, \frac{1}{2}]}: f(x) = |x|$.
So the inequality becomes an equality for all maximizers $f \in \mathbb{F}$. 
That implies that the only maximizers are functions that are linear with fixed slope $\pm 1$ in $[-\frac{1}{2}, 0] $ and $[0, \frac{1}{2}]$.
So the only possible optimizers are in the equivalence class $[|\cdot| |_{[-\frac{1}{2}, \frac{1}{2}]}]$ and $[id |_{[-\frac{1}{2}, \frac{1}{2}]}]$.
A short calculation reveals that the latter is not an optimizer.\\
Now define $X = \mathbb{F}, Y = P[-1/2, 1/2]$, $Y_1 = \{\pi_{1,\epsilon}|1 > \epsilon > 0 \}$ and $Y_2 = \{\pi_{2,\epsilon}|1 > \epsilon > 0 \}$.
Choose 
\[
\gamma(t) = \begin{cases}
\pi_{1,{\frac{1}{2}-t}},& t \in [0, \frac{1}{2}]\\
\pi_{2,t-\frac{1}{2}},& t \in [\frac{1}{2}, 1]\\
\end{cases},
\]
which is a computable path.
It holds
\[
\inf_{y_1 \in Opt(Y_1), y_2 \in Opt(Y_2)} \| y_1 - y_2 \|_{L^2([-\frac{1}{2}, \frac{1}{2}])} = \frac{\sqrt{5}}{4\sqrt{3}}.
\]
Note that $Y_1 \subset Y_1 \cup Y_2$ is decidable by checking if $f(1/2)$ is positive or negative. 
Now we can use Theorem \ref{mainthm}. 
This time we have to be careful since Theorem \ref{mainthm} only holds in the finite-dimensional case $X \subset \R_c^n$, $Y \subset \R_c^m$.
We can easily expand Theorem \ref{mainthm} for general Banach spaces with \textit{computable structure} using identical arguments.
For more information on computability on Banach spaces we refer to \cite{pour-el} and \cite{weihrauch}.
\end{proof}

\subsection{Lattice Problem for Cryptographic Applications}
\begin{proof}[Proof of Theorem \ref{thm:sivp}]
We slightly change the notation here.
Instead of bases, we consider \textit{ordered} bases.
This difference is only semantic and allows us to apply Theorem \ref{mainthm} more easily and does not change the nature of the optimization problem.
Define the bases $B_\lambda^1 = ((1,1), (-1-\lambda, 0) )$, $B_\lambda^2 = ((-1-\lambda, 0),(1,1) )$ for $\lambda \in [0,1/2]$.
Then for $\lambda = \sqrt{2}-1$ the minimizing bases are
\[
((-\lambda, 1), (1,1) ) \text{ and } ((-\lambda, 1), (1+ \lambda, 0) )
\]
as well as
\[
( (1,1), (-\lambda, 1)) \text{ and } ( (1+ \lambda, 0),(-\lambda, 1)).
\]
For $\lambda \in [0, \sqrt{2} - 1)$ the minimizing bases are
\[
((-\lambda, 1), (1+ \lambda, 0) ) \text{ and } ((1+ \lambda, 0), (-\lambda, 1) ).
\]
For $\lambda \in (\sqrt{2} - 1, 1/2]$ the minimizing bases are
\[
((-\lambda, 1), (1,1) ) \text{ and } ( (1,1), (-\lambda, 1) ).
\]
Now define the computable path $\gamma:[-1, 1] \to \{\text{ordered bases of }\R^2 \} $
\[
\gamma(x) := B_{x/100 + \sqrt{2} - 1 } \\
\]
So defining $X = Y = \{\text{ordered bases of }\mathbb{R}^2\} \subset \R^4$, $Y_1 = \{B^1_\lambda, B^2_\lambda  | \lambda \in (0, \sqrt{2} - 1) \}$ and $Y_2 = \{B^1_\lambda, B^2_\lambda  | \lambda \in (\sqrt{2}-1, 1/2) \}$, we see $\gamma([-1,0))\subset Y_1$ and $\gamma((0, 1])\subset Y_2$.
It holds
\[
\inf_{y_1 \in Opt(Y_1), y_2 \in Opt(Y_2)} \| y_1 - y_2 \|_{L^2([-\frac{1}{2}, \frac{1}{2}])} = \sqrt{2}.
\]
Now $Y_1 \subset Y_1 \cap Y_2$ is decidable by checking the first component of the vector $(-1-\lambda, 0)$ and comparing it with $\sqrt{2}$.
So by Theorem \ref{mainthm}, finding the optimizer is not Banach-Mazur computable.
 
\end{proof}

\subsection{Linear Program}
\begin{proof}[Proof of Theorem \ref{thm:lp}]
Define for $\epsilon \in (-1, 1)$ the linear Program
\begin{align*}
\textit{maximize  } & x_1 \\[1.5pt]
	\textit{subject to  }& x_1 \geq 0 \\
& x_1 + \epsilon x_2\leq 1 \\
& 0 \leq x_2 \leq 1 ,
\end{align*}
which can be coded by the parameter set 
\begin{align*}
c = &(1,0),\hspace{0.5cm} A_\epsilon = \begin{pmatrix}
-1 & 0 \\
1 & \epsilon \\
0 & -1\\
0 & 1
\end{pmatrix},\hspace{0.5cm} y = \begin{pmatrix}
0\\
1\\
0\\
1
\end{pmatrix}.
\end{align*}
We define $\Lambda_\epsilon := (c, A_\epsilon, y)$ as well as the disjoint sets $Y_i = \{\Lambda_{(-1)^{i+1}\epsilon} | \epsilon \in(0,1)\}$ for $i=1,2$.
Then $Opt(Y_1) = \{ (1,0) \}$ and $Opt(Y_2) = \{ (1 + \epsilon, 1)| \epsilon \in (0,1) \}$.
Define the computable path $\gamma: (-1,1) \rightarrow \mathbb{R}_c^n \times \mathbb{R}_c^{m\times n} \times \mathbb{R}_c^m$ by
\[
\gamma(t) = \Lambda_t.
\]
Now note that \\
\begin{enumerate}[label=(\arabic*)]
\item $Y_1 \subset Y_1 \cup Y_2$ is decidable by checking if $(A_\epsilon)_{1,1}$ is positive or negative.
\item $\gamma((-1,0)) \subset Y_2$ and $\gamma((0,1)) \subset Y_1$.
\item $\inf_{y_1\in Opt(Y_1), y_2 \in Opt(Y_2)} \|y_1 - y_2 \| = 1$.
\end{enumerate}
So by applying Theorem \ref{mainthm}, $G$ is neither Banach-Mazur computable nor $\alpha < \frac{1}{2}$ approximable by a Banach-Mazur computable function.
\end{proof}

\section{Acknowledgments}

Y. Lee acknowledges support by the German Research Foundation under Grant DFG-SPP-2298, KU 1446/32-1

This work of H. Boche and G. Kutyniok were supported in part by the ONE Munich Strategy Forum (LMU Munich, TU Munich, and the Bavarian Ministery for Science and Art).

This work of H. Boche was also supported in part by the German Federal Ministry of Education and Research (BMBF) in the project Hardware Platforms and Computing Models for Neuromorphic Computing (NeuroCM) under Grant 16ME0442 and within the national initiative on 6G Communication Systems through the research hub 6G-life under Grant 16KISK002.

G. Kutyniok acknowledges support from the Konrad Zuse School of Excellence in Reliable AI (DAAD), the Munich Center for Machine Learning (BMBF) as well as the German Research Foundation under Grants DFG-SPP-2298, KU 1446/31-1 and KU 1446/32-1 and under Grant DFG-SFB/TR 109, Project C09 and the Federal Ministry of Education and Research under Grant MaGriDo.


\bibliographystyle{unsrt}  
\bibliography{main}

\end{document}